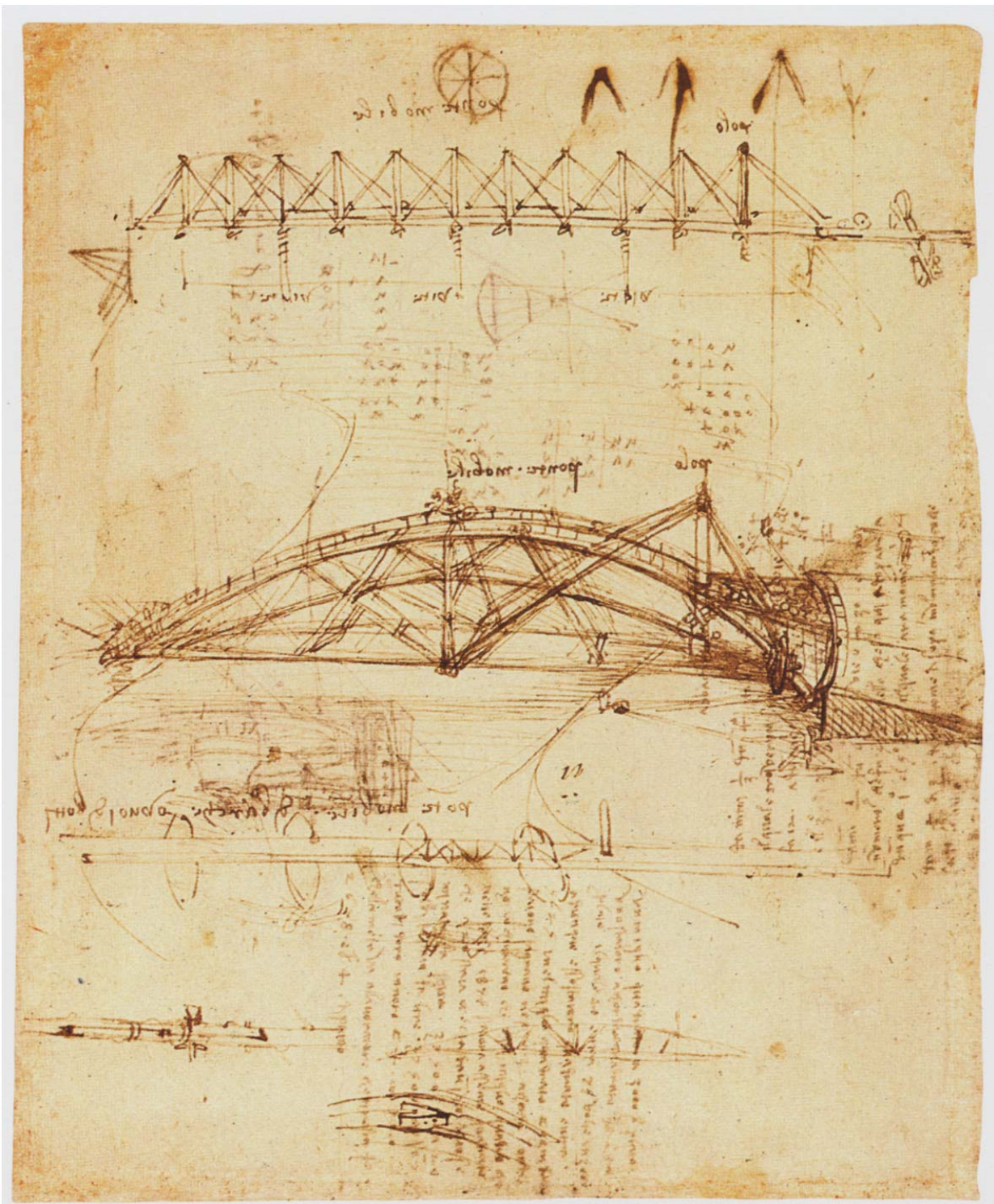

*Fig. 1: Folio 855 recto of Leonardo's Codex Atlanticus shows three movable bridges; the one on top is the problematic one.*

# A strange bridge by Leonardo

Dirk Huylebrouck,
LUCA/KULeuven, Paleizenstraat 65, 1030 Brussels, Belgium.

***This is an English translation of a paper that was submitted to the Dutch popular science magazine EOS (affiliated to Scientific American) in Sept 2012 and appeared in that magazine in April 2013.***

*On folio 855 recto of the Codex Atlanticus, Leonardo da Vinci drew three 'easily movable' bridges, but one of them is enigmatic: all 'replicas' in Leonardo museums and exhibitions come as a surprise, to say the least, to any engineer or architect whose attention is drawn to it. This is the case for models in Amboise (France), Chicago and Portland (USA), Florence (Italy) and for the one of the travelling exhibition by the Australian company 'Grande Exhibitions' that already visited 40 major cities in the world. All 'replicas' of the bridge model attributed to Leonardo have pillars standing on the deck of the bridge, while the deck is suspended by cables attached on these pillars. At first sight this problem does not catch the attention of the observer, as the bridge seems to be a mixed form of a beam and a suspension bridge, but it was not overlooked by my colleague architect-engineer Dr. Laurens Luyten (Gent, Belgium). Yet, after a TV-interview in Brussels so much pressure was exerted by some of the museum collaborators, the architect prefers to keep the silence about the engineering by these 'scientific' exhibitions. Yet, it is the author's opinion the ridiculous models do not honor Leonardo, and so the present paper takes up the challenge to report about his observation in his name (with his permission). Something should be done, for sure: millions of people look at the bridge replicas that immediately embarrass any engineer.*

**Introduction**

Three bridges are shown on folio 855 recto of the Codex Atlanticus, a 'box' with 1119 pages of drawings and writings by Leonardo da Vinci (Vinci, Italy, 1452 - Amboise, France, 1519), kept by his pupil Francesco Melzi. The pages were written between 1479 and 1519, but sculptor Pompeo Leoni later made some changes so that the current version of the codex in the Biblioteca Ambrosiana in Milan probably does not match the original compilation. The drawings on that folio show three movable bridges Leonardo proposed for the Borgia family (see fig. 1). He would have been asked to design light and strong bridges for their army so that soldiers could assemble and disassemble them quickly, using simple materials they could find or manufacture on the spot, such as logs and ropes. Thus, the goal was *not* necessarily to make bridges that could be moved along the shore, but bridges that could be easily built and rebuilt.

On the folio, Leonardo drew three bridges: one on eleven poles, a second supported by an arch, and a third resting on boats. The latter is a variant on the centuries-old principle of the pontoon bridge. The second construction method builds a bridge on one shore and swings the finished bridge deck around a single support standing on that shore until the other end rests on the opposite shore. This method recently enjoyed a revival for steel and concrete bridges, when a fast construction method is preferable due to heavy traffic on waterways.

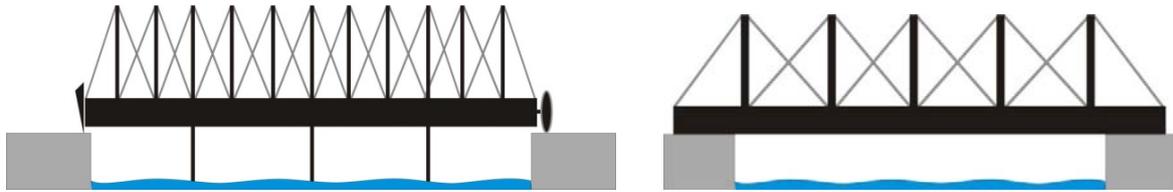

*Fig. 2: Schematic drawing of the first bridge by Leonardo on folio 855 recto of the Codex Atlanticus (left) and a common museum scale model of it (right).*

It is the first design that causes problems (see Fig. 2a), and more in particular how this bridge is shown in museum scale models, in books (see Fig. 2b, 3, 4 and 5) and, exceptionally, in a real-size model (see Fig. 6). Almost always, the number of pillars is not correctly copied. In addition, some models do not show any wheels at all, at none of the extremities. Even more surprising is their suspension method, using pillars standing on the deck of the bridge, while they are supposed to support the bridge itself.

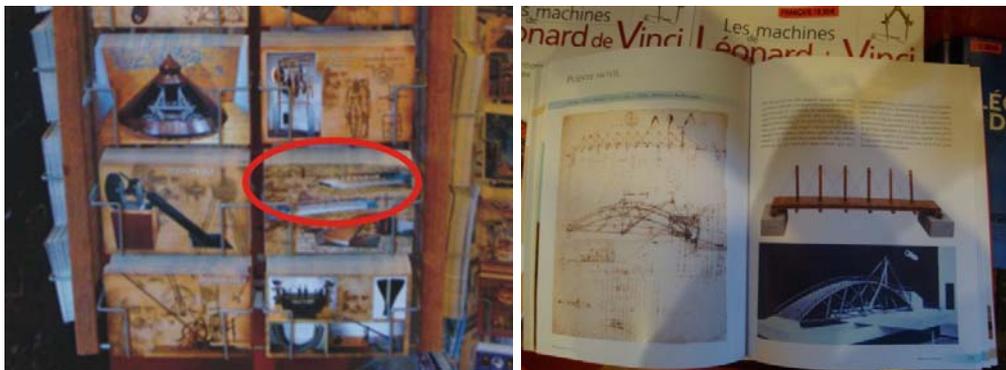

*Fig. 3: Shops in Leonardo museums, such as here in Amboise (France), sell postcards and books with an erroneous model for a movable Leonardo bridge.*

**A short-list of models**

The questionable model is for instance shown in the permanent exhibition 'Le Macchine di Leonardo Da Vinci' in Florence (Italy), the largest city near Vinci, the birthplace of Leonardo. And yet, Professor Carlo Pedretti, director of the Armand Hammer Center for 'Leonardo Studies' at the University of California at Los Angeles (USA), praised the accuracy of the models in this exhibition, though the example has only six pillars and wheels on both sides (see fig. 4). The exhibition room is currently undergoing a restoration and hopefully this paper will reach the organizers on time.

Almost identical models stand in the 'Museum of Science and Industry' in Chicago and in the 'Oregon Museum of Science and Industry', in Portland (USA), but in the latter case the bridge has only five pillars. An Australian company, 'Grande Exhibitions', rents travelling exhibitions to cities wanting to improve their tourist attractiveness. It has a Chinese Terracotta army exhibition, one about the Titanic, and another one showing Leonardo's work. The latter has the scientific support from the Italian organization 'Anthropos' and the Frenchman Pascal Cotte, and that is reasonable, regarding Leonardo's life. The exhibition currently remains in Brussels (Belgium), until September 2013, but later it will move to other towns and continents. In this exhibition, the scale model of the bridge we focus on here has only five pillars and no wheels at all (see fig. 4).

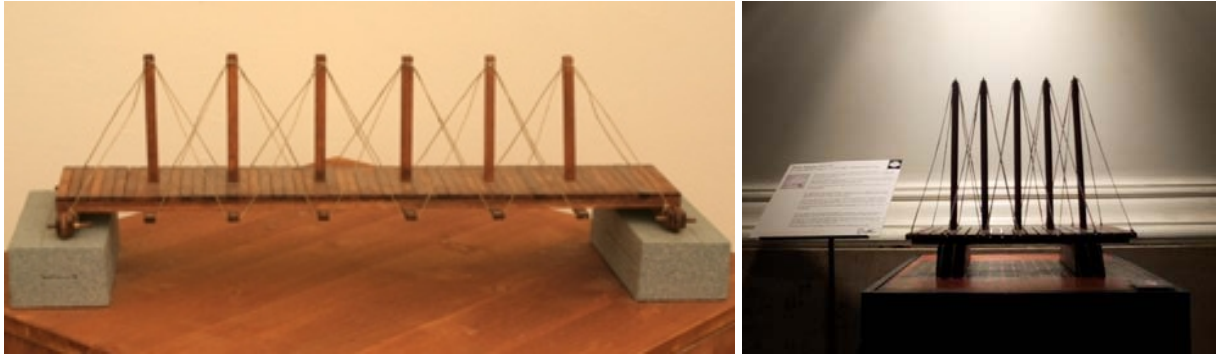
*Fig. 4: A six pillar 'model' with wheels at both sides as shown in Florence and in Chicago (left) and a five pillar model without any wheels shown by the travelling 'Grande Exhibitions', currently in Brussels (Belgium).*

In the 'Château du Clos Lucé' and its 'Parc Leonardo da Vinci', the number of columns, nine, comes the most close to that of Leonardo, eleven (see fig. 5). The 'Château' is located in Amboise, near Tours (France), the town where Leonardo died, three years after he was employed King Francis I. The 'mausoleum' shows a model with a small additional support below, in the middle, added by the model builder because there is no such support on Leonardo's sketch (the little support is not always in line with a pillar above as visitors touch the model). Perhaps the technician who made the model already experienced the problem of the bridge, even on a small scale: the bridge will inevitably bend if the deck is not very solid and so an additional support was necessary.

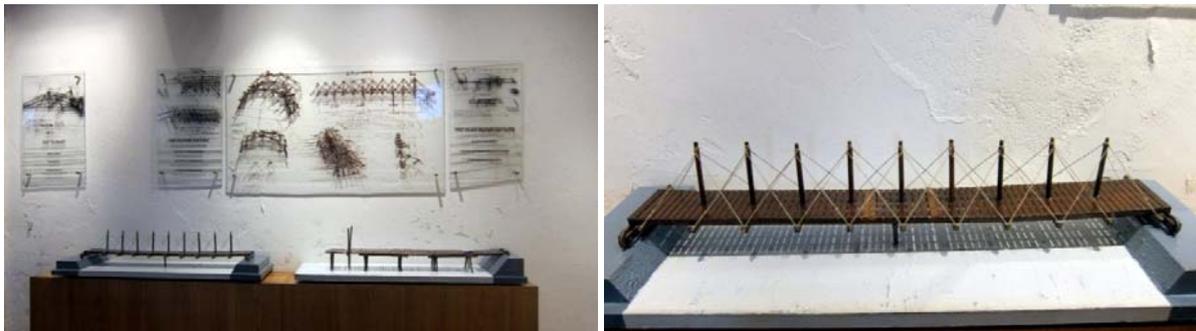
*Fig. 5: The presentation of the model in the museum of Amboise (left on the left picture) and a separate view (right).*

In the gardens surrounding the castle of Amboise, a large model of the bridge was build (see fig. 6). However, it has only two columns on each side and the wheels are purely decorative as they hang above the shore, and they are hard to spot. Also, the tops of the pillars are connected crosswise like in a truss bridge (see fig. 7 and 8) but there is not the slightest tension on the ropes as they are merely decorative. The explanatory plate with Leonardo's drawing stands about 100 meters away from the bridge at yet another bridge, different from the one it was supposed to explain (see fig. 8). Somebody messed up the panels, but this is understandable as the real size model of Leonardo's movable bridge does not correspond to Leonardo's drawing. The park designers must have been confused and put the panel with some other longer bridge.

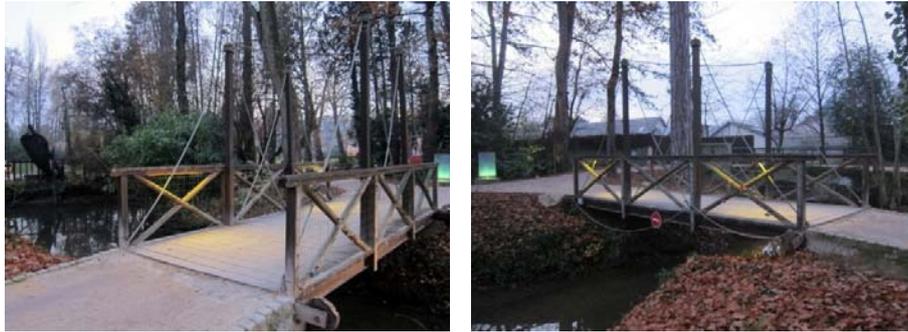

*Fig. 6: Two views on the 'true size' model in Amboise, with only 2 pillars on each side.*

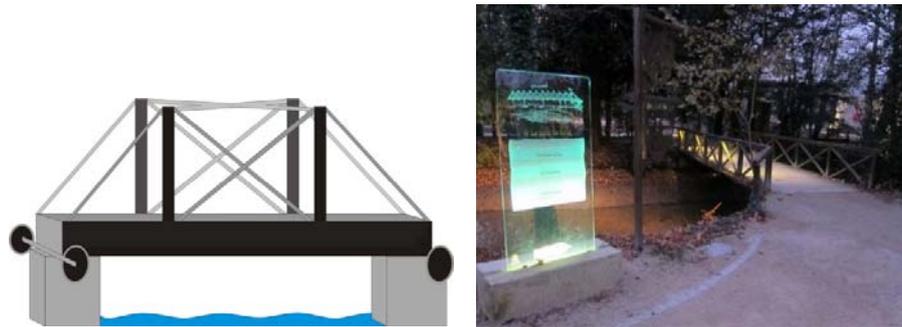

*Fig. 7: A schematic representation of the outdoor model in Amboise (left), but in the park the explanations accompany another bridge, some 100m away from the bridge it is supposed to illustrate (right).*

**The problem**

Leonardo's construction vaguely reminds the nowadays classical suspension bridge, and the similarity with such a hanging bridge and the scale models for Leonardo's movable bridge can be one of the reasons why the error does not immediately strike the museum visitor. Indeed, as any engineer or architect learns in a first year course on construction, there are several types of bridges. There is the 'beam bridge', which is as strong as the beam laying on the shores; the 'hanging bridge' where cables hanging on pillars allow a more slender deck; the 'cantilever bridge', using counterweights or fixations on the shores; the 'truss bridge' with cables connecting the tops of the pillars, thus considerably improving its strength; and the 'cable-stayed bridge', held up by cables like in a suspension bridge but requiring less cable for more explanations about these types of bridges (we omit the arc types here since there cannot be any confusion with those types). Of course, some of those types almost look like the scale models for Leonardo's bridge, and mixing one type with another can give the impression the scale models are correct.

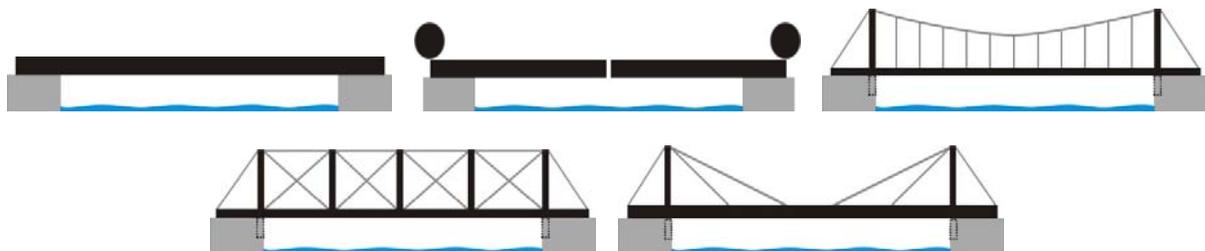

*Fig. 8: Five correct types of bridges: the beam bridge (above left), cantilever bridge (using spherical counterweights; above middle), suspension bridge (above right), truss bridge (below left) and cable-stayed bridge (below right).*

The common museum scale models for the Leonardo bridge incorrectly mix elements of these types of bridges: a strong deck as in a beam bridge, counterweights (the wheels) as in a cantilever bridge, and cables as in the suspension, truss or cable-stayed bridge. Yet, there always is an important difference: the ropes from the top of the pillars to the deck do not play any structural role in the scale models. The pillars all stand on the deck (and not on the shore as for the models in Fig. 8) and seem to lift themselves by ropes attached to that same deck. Thus, the scale models will not be stronger than the beam constituting the deck of the bridge. This can be checked using engineering software such as PowerFrame, though here its application is not really necessary (see fig. 9). Indeed, anybody with some common sense can immediately understand the problem and this probably explains why craftsmen made the scale models shorter and never with eleven pillars, though this was the number of pillars Leonardo drew.

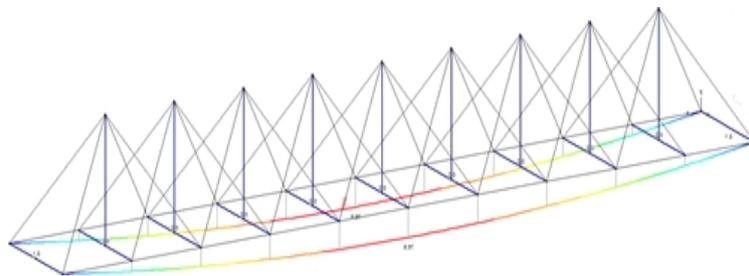

*Fig. 9: A computer study by L. Luyten of the Amboise model: there is engineering 'red alarm' in the deck, but there is no color change in the columns as they don't play a structural role.*

**Leonardo's idea?**
Perhaps the scale model designers were motivated to omit the parts standing in the river every three pillars because of the 'wheels' on the right side. They could suggest that the bridge had to be able to roll on the shore, but there are only 'wheels' on one side. The drawing is very clear and yet several models show wheels on both sides. Perhaps these 'wheels' were a kind of heavy millstones serving as counterweights at the start of the construction of the bridge, from the right to the left. It may have been the intention to build it beginning from the shore on the right of the drawing and move from the right to the left, advancing two 'modules' at the time in a cantilever construction. Next, from it a pillar would be pushed in the river bed and from this solid standing point, two more modules would again be built in cantilever in order to move to the other side, and so forth. In that case, the cables are pulled by the cantilever and that is indeed useful for the building process (see fig. 10). After all, an easy building process was Leonardo's goal, and not necessarily movability along the shore, as explained in the first section.

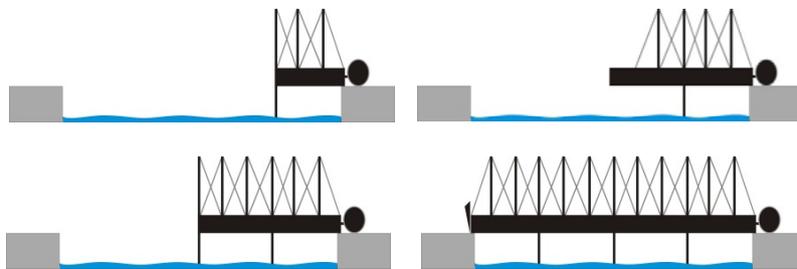

*Fig. 10: A possible cantilever construction, where the 'wheels' function as counterweights, like millstones. Progress is made from the right, pushing a pillar in the river from the cantilever every two pillars, as shown in Leonardo's drawing.*

Still, more guesswork remains. It difficult to explain why then the 'millstones' stand upright. The way in which they are attached suggests they are wheels and thus the interpretation of the museums to have no pillars standing in the river is understandable. In short, it remains unclear how Leonardo actually had in mind the first movable bridge on folio 855r. The construction method could also be imply an error was made by Leonardo, since he was not infallible in engineering, as Tarnai and Lengyel recently showed in other construction of his hand (see [3], [4], [5], [6]).

However, it should be noted his drawing is only a draft and not an official illustration in a book published by Leonardo. That applies to any drawing by Leonardo, because he never published a text or a book in his own name, and thus we will never know what his true intention was. In the Codex Arundel for instance figure strange erroneously drawn icosahedra, but of course Leonardo knew very well how to draw this polyhedron, since he did so for Luca Pacioli's '(De) Proportione Divina' (see Fig. 11, [2]).

Therefore, we will probably never have a precise idea about the first bridge designed on folio 855 of the Codex Atlanticus. But this should be no excuse for the numerous authors and museums to present obviously erroneous models. It can be acclaimed they try awaken the interest of the millions for mathematics and engineering, but they should respect the genius of Leonardo.

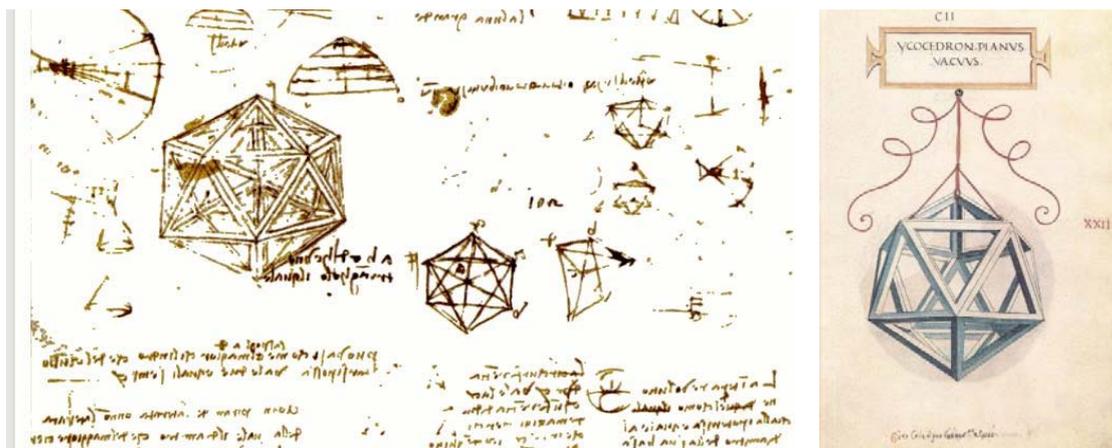

*Fig. 11: An erroneous draft drawn of an icosahedron (left) and a correct one (right), both attributed to Leonardo.*